\documentclass{amsart}
\usepackage{color}

\usepackage{amsmath} 
\usepackage{amssymb}
\usepackage{mathrsfs}

\newtheorem{theorem}{Theorem}[section] 
\newtheorem{claim}[theorem]{Claim}
\newtheorem{tic}[theorem]{The Isomorphism Claim}

\newtheorem{observation}[theorem]{Observation}

\theoremstyle{definition}
\newtheorem{definition}[theorem]{Definition}
\newtheorem{convention}[theorem]{Convention}

\newtheorem{fact}[theorem]{Fact}

\newtheorem{discussion}[theorem]{Discussion}

\newtheorem{hypothesis}[theorem]{Hypothesis}

\theoremstyle{remark}
\newtheorem{remark}[theorem]{Remark}

\newtheorem{notation}[theorem]{Notation}

\newcommand{\Sym}{{\rm Sym}}
\newcommand{\sym}{{\rm sym}}

\newcommand{\valency}{{\rm valency}}
\newcommand{\Prob}{{\rm Prob}}
\newcommand{\Graph}{{\rm Graph}}
\newcommand{\rang}{{\rm rang}}

\newcommand{\sg}{{\rm sq}}

\newcommand{\gr}{{\rm gr}}

\newcommand{\Rang}{{\rm Rang}}
\newcommand{\lqq}{{``}}
\newcommand{\snn}{{\sn}}
\newcommand{\EXP}{{\rm EXP}}

\newcommand{\rest}{{\restriction}}
\newcommand{\dom}{{\rm dom}}

\newcommand{\iif}{{\rm if}}

\newcommand{\then}{{\underline{then}}}
\newcommand{\when}{{\underline{when}}}
\newcommand{\oor}{{\underline{or}}}

\newcommand{\Iff}{{\underline{iff}}}
\newcommand{\but}{{\underline{but}}}

\newcommand{\mn}{{\medskip\noindent}}
\newcommand{\sn}{{\smallskip\noindent}}

\newcommand{\cE}{{\mathscr E}}
\newcommand{\cL}{{\mathscr L}}
\newcommand{\bbB}{{\mathbb B}}

\newcommand{\bbR}{{\mathbb R}}

\newcommand{\varp}{{\varepsilon}}

\newcommand{\cG}{{\mathscr G}}

\newcommand{\bbL}{{\mathbb L}}
\newcommand{\bbN}{{\mathbb N}}
\newcommand{\cM}{{\mathscr M}}

\newcommand{\bbQ}{{\mathbb Q}}

\newcommand{\cT}{{\mathscr T}}

\newcount\skewfactor
\def\mathunderaccent#1#2 {\let\theaccent#1\skewfactor#2
\mathpalette\putaccentunder}
\def\putaccentunder#1#2{\oalign{$#1#2$\crcr\hidewidth
\vbox to.2ex{\hbox{$#1\skew\skewfactor\theaccent{}$}\vss}\hidewidth}}

\newenvironment{PROOF}[2][\proofname.]
   {\begin{proof}[#1]\renewcommand{\qedsymbol}{\textsquare$_{\rm #2}$}}
   {\end{proof}}

\usepackage[hidelinks]{hyperref}

\begin{document}
\makeatletter\def\shfiuwefootnote{\gdef\@thefnmark{}\@footnotetext}\makeatother\shfiuwefootnote{Version 2021-07-15. See \url{https://shelah.logic.at/papers/1077/} for possible updates.}

\title {Random graph: stronger logic but with the zero one law \\
Sh1077}
\author {Saharon Shelah}
\address{Einstein Institute of Mathematics\\
Edmond J. Safra Campus, Givat Ram\\
The Hebrew University of Jerusalem\\
Jerusalem, 91904, Israel\\
 and \\
 Department of Mathematics\\
 Hill Center - Busch Campus \\ 
 Rutgers, The State University of New Jersey \\
 110 Frelinghuysen Road \\
 Piscataway, NJ 08854-8019 USA}
\email{shelah@math.huji.ac.il}
\urladdr{http://shelah.logic.at}
\thanks{The author thanks Alice Leonhardt for the beautiful typing.}





\date{2021-05-29a}

\begin{abstract}
We like to find a logic really stronger than first order for the
random graph with edge probability $\frac 12$ but satisfies the 0-1
law.  This means that on the
one hand it satisfies the 0-1 law, e.g. for the random graph
$\cG_{n,1/2}$ and on the other hand there is a 
formula $\varphi(x)$ such that for no first order $\psi(x)$ do we have:
for every random enough $\cG_{n,1/2}$ the formulas
$\varphi(x),\psi(x)$ 
are 
equivalent in it.  We do it adding a quantifier
on graphs $\bbQ_{\mathbf t}$, i.e. have a class of finite graphs closed
under isomorphisms and being able to say that if $(\varphi_0(x,\bar
c),\varphi_1(x_0,x_1,\bar c))$, a pair of formulas with parameters
defining a graph in $\cG_{n,1/2}$, \then \, we can form a
formula $\psi(\bar y)$ such $\psi(\bar c)$ 
says that the graph belongs $K_{\bar{\mathbf t}}$.  
Presently we do it for random enough $\bar{\mathbf t}$.  
\smallskip

\noindent



\end{abstract}

\maketitle
\numberwithin{equation}{section}
\setcounter{section}{-1}
\newpage

\centerline {Annotated Content}
\bigskip

\noindent
\S0 \quad Introduction, pg.\pageref{0}
\bigskip

\noindent
\S1 \quad Identifying the too simple graphs, pg.\pageref{1}
\mn
\begin{enumerate}  
\item[${{}}$]  [We choose a $\mathbf h:\bbN \rightarrow (0,1)_{\bbR}$
 going to zero slowly enough.  Out intention is to add to first-order
logic a quantifier describing random properties of a graph but excluding
 some ``low", ``explicitly not random" graphs.
 Those are 
graphs such that for any quantifier free first order formula $\varphi(\bar
 x_0,\bar x_1,\bar z)$ for some $k$, for random enough $G =
 \cG_{n,1/2}$ (or $\cG_{n,p}$ for a given $p \in (0,1)_{\bbR}$), 
if $\bar c \in {}^{\ell g(\bar z)}G$ and
$\varphi(\bar x_0,\bar x_1,\bar c)$ define in $G$ a graph with $> k$
 nodes then it is so called low.  This will be used in \S2 to find a logic
 as desired.]
\end{enumerate}
\bigskip

\noindent
\S2 \quad The Quantifier, pg.\pageref{2}
\mn
\begin{enumerate}  
\item[${{}}$]  [We choose randomly enough a set $\mathbf K$ of (isomorphism
types of) finite non-$\mathbf k$-low graphs and show that adding a
  quantifier for it preserves the zero-one law. 
  So,  
  the probability of $H$, a non-low graph to be 
in the class is $\mathbf h(|H|)$.
  Why $\mathbf h$ is not constant?  
Because we like that on the one hand, $\Pr(\cG_{n,p} \in \mathbf K)$ 
converge to 0 (or to 1) so that a sentence saying 
(the graph $\cG_{n,p}$ belongs to $\mathbf K$)
converge to 0 (or to 1), and 
for   
similarly any graph definable in
$\cG_{n,p}$ by a first order formula without parameters.  On the other
hand, the probability of e.g. ``there is $a \in \cG_{n,p}$ such that
$\cG_{n,p} \rest \{b:bRa\}$ belongs to $K$"  
will go to 1, as there
will be $\ll n$ such nodes but still many.]
\end{enumerate}
\newpage

\section {Introduction} \label{0}

Our aim is to find a logic $\cL$ stronger than first order such that:
for $p \in (0,1)_{\bbR}$, the $p$-random graph $\cG = \cG_{n,p}$ (i.e. with edge
probability $p$) satisfies the 0-1-law
\underline{but} some formula $\varphi(x) \in \bbL$(graphs) 
defines in random enough graph $\cG_{n,p}$ a set of 
nodes not definable by any first order logic formula (of course, small enough
compared to $n$, even with parameters).

The logic is gotten from first order $\bbL$ by adding a (Lindstr\"om)
quantifier $\mathbf Q_{\bar{\mathbf t}} = \mathbf Q_{\mathbf{K} _{\bar{\mathbf t}}}$
gotten from a ``random enough" $\bar{\mathbf t} \in 
{}^{\bbN}\{0,1\}$; on quantifiers see \cite{BF}.  
We may wonder, can we replace $\mathbf Q$ by a ``reasonably defined
quantifier"?  We may from the proof see what we need from $\mathbf K$,
the class defining the quantifier $\mathbf Q_{\mathbf K}$, i.e.
a class of (finite) graphs closed under isomorphisms.
Excluding some graphs which we call low, the
membership in $\mathbf K$ should 
be random enough in the sense that
if we consider only random enough $\cG_{n,p}$, the non-trivial
$\bbL(\mathbf Q_{\mathbf K})$-formulas with parameters 
will define graphs which are not low and
are pairwise non-isomorphic except in trivial cases.  So we just need a
definition satisfying this; we hope  
to try to do it in a work in preparation.

How does the randomness of $\bar{\mathbf t}$ help us to 
get the zero-one law?  The idea is that for 
the quantifier $\mathbf Q_{\bar t}$ (see \S2) used here, if we
expand $\cG_{n,p}$ by finitely many relations definable by
formulas from $\bbL(\mathbf Q_{\bar{\mathbf t}})$, we get a random
structure with more relations essentially with constant probabilities, i.e. is
interpretable in a  
suitable $\cM =  
\cM_{\mathbf s,\bar p,n}$, see \S1,  
it look like $ {\mathscr G} _ 
        {p,n}$ (but with some 
relations of suitable kinds as we sort out), 
with, e.g. $\bar p 
= \langle p_n:n < \omega\rangle$ with $p_n$ going slowly to zero.

That is, fixing formulas $\varphi_\ell(\bar x_\ell) \in \bbL(\mathbf
Q_{\bar{\mathbf t}})$ starting with $\mathbf Q_{\bar{\mathbf t}},\ell
< k$ with no obvious connections we decide a priory that for a
random enough 
$\cG_{n,\bar p}  $ 
the structure $  
    ([n],R^{\cG_{n,p}}_\ell)_{\ell < k}  
    =( \varphi ({\mathscr G} _{p.n}), \dots ,   
        \varphi _ {\ell} ({\mathscr G} _{p,n}), \dots )_ {\ell} $   
        for suitable formulas 
        $ \varphi (x), \varphi _{\ell} 
    ( \bar{ x } _ {\ell} )$, 
        will look like $ {\mathscr M } $ above.

The decision is the simplest one: look as if
truth values of $R^{\cG_{n,p}}_\ell(\bar a)$ were drawn independently,
with probability $p_n$.  This is an over simplification!  We need a
more involved such drawing, reflecting the original $\bar\varphi_\ell$
to some extent, see below.

We may replace $\cM_{\mathbf s,\bar p,n}$ by using 
(for some irrational $\alpha \in
(0,1))$ $\bar p_n = (p,p_n)$, such that $p_n =
\frac{1}{n^\alpha}$, except the original drawing of the graphs as in
\cite{Sh:304}.  We can also analyze 
    $ {\mathscr G} {n, r n^ \alpha }$  
and use
several pairs $(r,\alpha)$ 
in the analysis (as long as the sets of
$\alpha$'s is linearly independent over the rationals).  We 
hope  
later to show that for
some such version there is a more natural definable $\mathbf Q_{\mathbf K}$
which imitate its behavior.

So in the proof we have two questions to address: first fixing $G =
([n],R_\ell)_{\ell < k}$, 
drawing the quantifiers, how $([n]   
    ,R^G_\ell,\ldots)$ look
like.  Second, we need to consider all the $G$'s on $[n]$.  For the
first stage the main problems are: two definably derived graphs which 
are isomorphic.

We do some kind of elimination of quantifiers: essentially if $\cM_n$ is
a $\tau$-structure ($\tau$ relational and finite) drawn randomly
according to the sequence $\langle p_{\tau,R}:R \in \tau\rangle$ of
fixed probabilities, applying $\mathbf Q_{\bar{\mathbf t}}$ to some finitely many
schemes $\langle \mathbf s_0,\dotsc,\mathbf s_k\rangle$ of interpreting
graphs, define a random $\cM'_n$ for $\tau'$-structures by
expanding $\cM_n$ by $R_\ell = \{\bar c:\ell g(\bar c) = \ell g(\bar
z_{\mathbf s_\ell})$ and the graph $H_{\mathbf s_\ell,\bar c}$ interpreted
by $\mathbf s_\ell$ for the parameter $\bar c$ is in the class 
$\mathbf Q_{\bar{\mathbf t}}\}$.  

Our use of vocabulary and structure deviates a little from the
standard, but fits with the use in graph theory and is natural
here.  In graph theory the edge relation $R$ is assume to be symmetric
and irrefelxive.  So we use (say $k_t$-place predicate) $R_t$ such that
it is always irreflexive (fails for $k_t$-tuples with a repetition)
and $K_t$-invariant for some group $K_t$ of permutation of
$\{0,\dotsc,k_t-1\}$, i.e. if $\langle a_\ell:\ell < k_t\rangle$
satisfies it then so does $\langle \bar a_{\pi(\ell)}:\ell <
k_t\rangle$ for every $\pi \in K_t$.  This is natural because when the
pair of formulas 
$\bar\varphi(\bar c)$ defines a graph $H=H_{M,\bar\varphi,\bar c}$ in
the structure $M$ (e.g. a graph) and we like to draw a truth value for
``$H \in \mathbf K_{\bar{\mathbf t}}$", a group of permutation of $\ell
g(\bar c)$ is dictated by $\bar\varphi$.

Why the random auxiliary structures are better defined in a different
way?  Recall the truth value of ``$H \in \mathbf K_{\bar{\mathbf t}}$" is
chosen randomly, but if $H$ is definable in the graph $G$, say is
$H_{G,\bar\varphi,\bar c,\bar{\mathbf t}}$ then the probability of ``$H \in \mathbf
K_{\mathbf t}$" depends on $H$, and in natural cases, on $|H|$, the
number of nodes of $H$. But if $\cM = ([n],\dotsc,R^0_\ell,\ldots)$ is
random, the standard way  
to make the probability of $\bar c \in
R^G_\ell$ naturally depend on $n$ and in many cases $n \ne |H|$.

We could have allowed using the quantifiers only on graphs $H$
definable in $\cG_{n,q}$ with set of nodes $[n]$ but this seems to me
quite undesirable, 
restricting our logic too much.     
We restrict ourselves to the class of graphs -
twice, we consider $\cG_{n,q}$ and the quantifier $\mathbf Q_{\bar{\mathbf
t}}$ is on graphs.  But in both cases this is not really needed.

We thank Simi Haber for raising again the problem and for some
stimulating discussions and Noga Alon for asking during a lecture in
the Noga-fest, January 2011, why we ignore the weak graph; 
a reasonable interpretation is:  
why we do not draw a truth value for \lqq $ G $ is green" 
for $ G $ a empty  graph. One problem is that the sentence $ \psi $ 
saying  \lqq the  graph  with all nodes (is $[n]$) and no edges"
the  
probability that $ {\mathscr G} _{p, n }$ satisfies
it is always zero or one and in non-trial cases
is not eventually constant; see more in \S3

\newpage

\section {Identifying the low graphs} \label{1}

We like to add a quantifier $\mathbf Q$ on finite graphs, which give a
property of finite graphs respecting isomorphism
(i.e.  a subset closed under automormorphisms).  
The aim is that for e.g. for the random graph $\cG_{n,p}$, 
the 0-1 law holds for $\bbL(\mathbf Q)$ but there is an $\bbL(\mathbf
Q)$-formula $\varphi(x)$ such that for no first order $\psi(x)$ are
$\varphi(x),\psi(x)$ equivalently in $\cG_{n,p}$.

More specifically, we
better make the quantifier trivial on too simple graphs, then we
intend that for any
fix finite set of formulas from $\bbL(\mathbf Q)$, for random enough
$G_{n,p}$ the structure $(G,\varphi^G(-))_{\varphi \in \Delta}$ is a
random structure excluding the ``problematic" graphs.
\bigskip

\subsection {Interpretation}\
\bigskip

\begin{convention}
\label{a0}
1) $\mathbf h:\bbN \rightarrow (0,\frac 12)_{\bbR}$ goes to zero slowly
enough, e.g. $\mathbf h(n) = 1/\log_2 \log_2(n)$ for $n > 16$ and $=1$ if
$n \le 16$; slowly enough actually means:
\mn
\begin{enumerate}
\item[$(a)$]    
$\alpha \in (0,1)_{\bbR} \Rightarrow 
\infty   
= \lim\langle \mathbf h(n) n^ \alpha:n < \omega\rangle$  
\sn
\item[$(b)$] $0=  \lim\langle \mathbf h(n):n < \omega\rangle$  
\snn 
\item[$(c)$] $ n^{\mathbf{h}(n) }$   
is non-decreasing (for simplicity).   
\end{enumerate}
\mn
2) $\mathbf g:\bbN \rightarrow \bbR_{\ge 0}$  
be $\mathbf g(n) = 
n^{\mathbf h(n)}$ hence $\mathbf g(1+n) \ge 1$ 
and so 
$\mathbf g$ go to 
infinity slowly enough.
\end{convention}

\begin{notation}
\label{a1}
1) Let $[n] = \{1,\dotsc,n\}$ or $\{0,\dotsc,n-1\}$ if you prefer (serve
as the universe of the $n$-th random graph).

\noindent
2) Let $\tau$ denote a vocabulary (e.g. $\tau = \tau_{\gr}$ is
the vocabulary of graphs; see Definition \ref{a2} below).  Let $\bbL$
be first order logic so $\bbL(\tau)$ is the set of first order
formulas in the vocabulary $\tau$, but below we may write $\bbL(\mathbf
s)$ instead of $\bbL(\tau_{\mathbf s})$.

\noindent
3) A $\tau$-model $M$ is defined as usual.

\noindent
4) For a formula $\varphi = \varphi(\bar x,\bar y)$, model $M$ and
$\bar b \in {}^{\ell g(\bar y)}M$ let $\varphi(M,\bar b) = \{\bar a
\in {}^{\ell g(\bar x)}M:M \models \varphi[\bar a,\bar b]\}$.
\end{notation}

The following is a central definition, explicating the restriction 
    to what is definable.  
\begin{definition}  
\label{a2}
1) For a finite set $I$ we say $\mathbf s$ is an $I$-kind or an $I$-kind
sequence (of a vocabulary) and write $I_{\mathbf s} = I$ \when \,:
\mn
\begin{enumerate}
\item[$(a)$]  $\mathbf s= \langle(k_t,K_t):t \in I\rangle = \langle
  (k_{\mathbf s,t},K_{\mathbf s,t}):t \in I\rangle$
\sn
\item[$(b)$]  $k_t \in \bbN$
\sn
\item[$(c)$]  $K_t$ is a group of permutations of
  $\{0,\dotsc,k_t-1\}$.
\end{enumerate}
\mn
1A) Let $\mathbf s_{\gr} = \mathbf s(\gr)$ be defined by ($\gr$ stands for
graphs) $I_{\mathbf s} = \{s_{\gr}\},s_{\gr}$ fix, e.g. $0,k_{\mathbf
  s,s_0} = 2,K_{\mathbf s,s_0} = \sym(2)$, the group of permutations of
$\{0,1\}$.

\noindent
2) For $\mathbf s$ an $I$-kind sequence we define:
\mn
\begin{enumerate}
\item[$(a)$]  the $\mathbf s$-vocabulary $\tau_{\mathbf s}$ is $\{R_t:t
  \in I\},R_t$ a $k_{\mathbf s,t}$-place predicate
\sn
\item[$(b)$]  an $\mathbf s$-structure is $M=(|M|,R^M_t)_{t \in I}$ such
  that (so the universe $|M|$ of $M$ may be empty):
\sn
\begin{enumerate}
\item[$(\alpha)$]  $R^M_t$ is a $k_t$-place relation on $|M|$
\sn
\item[$(\beta)$]  $R^M_t$ is $K_t$-invariant, i.e. if $\langle
  a_\ell:\ell < k_t\rangle \in R^M_t \wedge \bar b \in 
\bar a/E_{K_t} \Rightarrow \bar b \in R^M_t$ where $\bar
  a/E_{K_t} = \{\langle a_{\pi(\ell)}:\ell < n_t\rangle:\pi \in
K_t\}$; let $E_{\mathbf s,t} = E_{K_t}$
\sn
\item[$(\gamma)$]  $R^M_t$ is irreflexive, i.e. $\bar a \in R^M_t
  \Rightarrow \bar a$ with no repetitions.
\end{enumerate}
\sn
\item[$(c)$]  $\mathbf M_{\mathbf s} = \cup\{\mathbf M_{s,n}:n \in \bbN\}$
  where $\mathbf M_{s,m} = \{M:M$ an $\mathbf s$-structure with set of
  elements $[n]\}$.
\end{enumerate}
\mn
3) For an $I$-kind $\mathbf s$ let $\mathbf P^1_{\mathbf s}$ be the set of 
$\bar p = \langle p_{t,n}:t \in I,n \in \bbN\rangle,p_{t,n} 
\in (0,1)_{\bbR}$, so $p_{t,n} \notin \{0,1\}$.  
We define the $(\mathbf s,\bar p)$-random
structure on $[n],\cM = \cM_{\mathbf s,\bar p,n}$
as follows (see more in part (5),(6):  
for $t \in I$ and $\bar a \in {}^{k_t}([n])$ with 
no repetitions we draw a truth value for $\bar a
\in R^{\cM}_t$ with probability $p_{t,n}$, \underline{but}
demanding we have the same result for $\bar a',\bar a''$ when 
they are $E_{\mathbf s,t}$-equivalent and independent otherwise.

\noindent
3A) Let $\mathbf P^0_{\mathbf s}$ for $\mathbf s$ as above be the set of
$\bar p \in \mathbf P^1_{\mathbf s}$ such that $t \in I_{\mathbf s} \wedge n \in \bbN
\Rightarrow p_{t,n} = p_{t,0}$, so we may write $p_t$ instead of
$p_{t,0}$.  If $\mathbf s = \mathbf s_{\gr}$, we may 
write $ \gr $ instead $ \mathbf{s} $.  

\noindent
4) Let $\mathbf P^2_{\mathbf s}$ be the set of $\bar p \in \mathbf
P^1_{\mathbf s}$ such that for some $\bar q \in \mathbf P^0_{\mathbf s}$ and
partition $\bar I = (I_0,I_1)$ of $I$, we have
$p_{t,n}$ is  
$q_t$ if $t  
\in I_0$ and is $q_0/\mathbf g(n)$ if $t \in
I_1$; we denote $\bar p$ by $\bar p_{\bar q,I_0} = \bar p_{\bar q,\bar I  
} = \bar p[\bar q,\bar I]$.

\noindent
4A) We may write $p$ instead of $\langle p_t:t \in I_{\gr}\rangle$
when $p_{s_{\gr}} =p$.

\noindent
5) For $\bar p \in \mathbf{P}^1_ \mathbf{s} $ let 
 $\mu_{\mathbf s,\bar p,n}$ be the distribution (= probability
   space) on $\mathbf M_{\mathbf s,n}$ corresponding to drawing the truth
   value of $R_t(\bar a)$ really of $\langle \mathbf R_t(\bar a'):
\bar a' \in \bar a/E_{\mathbf s,t}\rangle$ for a sequence $\bar a$ 
with no repetitions of length $k_{\mathbf s,t}$ with probability
   $p_{t,n}$, independently of the other choices.

\noindent
6) Let $\cM_{\mathbf s,\bar p,n}$ be the random variable for
the finite probability space $(\mathbf M_{\mathbf s,n},\mu_{\mathbf s,\bar p,n})$.
If $\mathbf s = \mathbf s_{\gr,q}$ 
let $ {\mathscr G} _  
    {q,n} = 
\cG_{\gr,q,n} =  
\cM_{\mathbf s_{\gr},q,n}$ and $\mu_{\gr,q,n} = 
\mu_{\mathbf s_{\gr},\bar p_{\gr,q},n}$.
\end{definition}

\noindent
Recall
\begin{fact}  
\label{a4}
1) $\mathbf P^0_{\mathbf s} \subseteq \mathbf P^2_{\mathbf s} \subseteq \mathbf
P^1_{\mathbf s}$.

\noindent
2) For every $\bar p \in \mathbf P^0_{\mathbf s}$ or $\bar p
 \in \mathbf P^2_{\mathbf s},\cM_{\mathbf s,\bar p,n}$ 
satisfies the 0-1 law for first order logic
and the limit theory $T_{\mathbf s,\bar p}$ has elimination of
quantifiers, really is $T_{\mathbf s}$, i.e. does not depend on $\bar p$
and $\mathbf g$ and $\mathbf h$ (as long as they are as in \ref{a0}(2)).

\noindent
3) $\mathbf M_{\mathbf s_{\gr},n}$ is the set of graphs with set of nodes
   $[n]$.
\end{fact}

\begin{PROOF}{\ref{a4}}
Should be clear.
\end{PROOF}

\begin{remark}  
\label{a5} 
1)
We first concentrate on one application of the quantifier.

\noindent
2)
We are interested in interpreting graphs.  We give the most general
case.  Note that 
we intend the quantifier to be a property of
  graphs.  So we have to think of an interpretation of a graph.  In
  such general interpretations using quantifier free formulas the
 elements may be only: a set of elements definable by a formula
$\varphi(x,\bar a),\bar a$ is a sequence of parameters \underline{or} more
generally such a set of $k$-tuples, maybe modulo suitable $E_K$,
  or even a finite union of such.  For each pair of the nodes
(fixing from where in the union they come) we define when it is an
edge by a quantifier free formula.
So below $\bar z$ are parameters, $\mathbf i(\bar\varphi)$ number of ``kinds
of elements",  
ways to define a node; $\varphi_{0,i}$ restrict 
the $i$-th kind,
$\varphi_2(\bar z)$ describes the relevant parameters, $\varphi_{1,i,j}$
describes the edges between a node of the $i$-th kind and a node of
the $j$-th kind.

\noindent 
3) Generally in interpretations we allow the set of elements
to be e.g. the set of equivalence classes of 
an equivalence relation 
defined say by $ \varphi (\bar{ x } ', \bar{ x } '', \bar{ a } )$, 
    where $\lg( \bar{  x }' ) = \lg( \bar{ x}'' )$
\underline{but} in our case those will always be degenerated,  
see \ref{a15}.  
\end{remark} 

\begin{definition}  
\label{a6}
1) For $\mathbf s$ an $I$-kind, we say 
$\bar\varphi$ is a $\mathbf s$-scheme (of a graph interpretation in
$\mathbf s$-structures) \when \, it consists of:
\mn
\begin{enumerate}
\item[$(a)$]  $\langle \varphi_{0,i}(\bar x_i,\bar z),\varphi_{1,i,j}(\bar
  x_i,\bar x'_j,\bar z),\varphi_2(\bar z):i,j < \mathbf
  i(\bar\varphi)\rangle$ such that:
\sn
\begin{enumerate}
\item[$\bullet_1$]   $\ell g(\bar x'_j) = \ell g(\bar x_j)$, it is
  possibly zero 
\sn
\item[$\bullet_2$]   $\langle \bar x_i,\bar x'_i:i < \mathbf
  i(\varphi)\rangle$ are pairwise disjoint, each with no repetitions
\sn
\item[$\bullet_3$]   $\mathbf i(\bar\varphi)$ is a non-zero natural number;
    if we allow $ \mathbf{i} ( \bar{ \varphi } ) =0 $ then we have to allow
        the empty graph.  
\end{enumerate}
\sn
\item[$(b)$]  $\varphi_{0,i},\varphi_{1,i,j},\varphi_2$ are formulas
  in the vocabulary $\tau_{\mathbf s}$, in this section they 
  always are  
 quantifier free formulas in $\bbL(\tau_{\mathbf s})$, 
    note that possibly $ \varphi _{1,i= \varphi_{1,j}}$ though  
    $i\not=  j $. 
\sn
\item[$(c)$]  $K_i = K_{\bar\varphi,i}$ is a group of 
permutations of $\{0,\dotsc,\ell g(\bar x_i)-1\}$, 
not related to $K_{\mathbf s,t}(t \in I)!$
\sn
\item[$(d)$]  $\varphi_{0,i}(\bar x_i,\bar z)$ is invariant under
permuting $\bar x_i$ by any $\pi \in K_i$;
    that is 
if $\pi \in K_i$; $\bar
x'_i = \langle x_{i,\pi(f)}:\ell < \ell g(\bar x_i)\rangle$ then
$\varphi_2(\bar z) \vdash_{\mathbf s} (\forall x_0 \ldots x_\ell
\ldots)(\varphi_{1,i}(\ldots x_\ell \ldots;\bar z) \equiv
\varphi_{1,i}(\ldots x_{\pi(\ell)},\dotsc,\bar z))$ where
$\vdash_{\mathbf s}$ means implication in every $\mathbf s$-structure
\sn
\item[$(e)$]  $\varphi_{1,i,j}(\bar x_i,\bar x'_j,\bar z)$ is
  invariant under permuting $\bar x_i,\bar x'_j$ by $\pi \in
  K_i,\varkappa \in K_j$ respectively, and 
$\vdash \varphi_{1,i,j}(\bar x_i,\bar x'_j,\bar z) 
\equiv \varphi_{1,j,i}(\bar x'_j,\bar x_i,\bar z)$ and $\vdash \neg
\varphi_{1,i,i}(\bar x_i,\bar x_i,\bar z)$
\sn
\item[$(f)$]  if $M$ is a 
$\tau_{\mathbf s}$-structure and $G \models \varphi_{0,i}[\bar a,\bar
c]$, so $\ell g(\bar c) = \ell g(\bar z)$ \then \, $\bar a 
\char 94 \bar c$ is with no repetitions.
\end{enumerate}
\mn
So if we have $\bar\varphi = \bar \varphi^\iota$ then
$\varphi^\iota_{0,i} = \varphi_{0,i}$, etc. and we may write $\bar
z_{\bar\varphi},\bar x_{\bar\varphi,1,i},\bar x'_{\bar\varphi 
,1,1}$.

\noindent
2) If $\mathbf s$ and $\bar\varphi$ are as above, $M$ is an $\mathbf
s$-structure and 
$\bar c \in \varphi_2(M)$, i.e. $\bar c \in {}^{\ell g(\bar z)}M$ 
satisfies $M \models \varphi_2[\bar c]$ \then \, 
$H=H_{\bar\varphi,\cM,\bar c}$ is 
    the follwoing 
    graph:  
\mn
\begin{enumerate}
\item[$(\alpha)$]  the set of nodes is $\{(i,\bar a/E_{K_i}): M \models
\varphi_{0,i}[\bar a,\bar c]$ for some $i < \mathbf i(\bar\varphi)$ and 
$\bar a \in {}^{\ell g(\bar x_i)}M\}$, see \ref{a2}(2)(b)$(\beta)$
\sn
\item[$(\beta)$]  $\{(i,\bar a/E_{K_i}),(j,\bar b/E_{K_j})\}$ is an edge 
iff $M \models \varphi_{1,i,j}[\bar a,\bar b,\bar c]$.
\end{enumerate}
\mn
3) Let $k_{\mathbf s}(\bar\varphi) = \max(\{\ell g(\bar x_i):i < \mathbf
i(\bar\varphi)\} \cup \{\ell g(\bar z)\})$ and let $k_{\mathbf
  s,i}(\bar\varphi) = \ell g(\bar x_i),k^*_{\mathbf s}(\bar\varphi) =
\max\{\ell g(\bar x_i):i < \mathbf i(\bar\varphi)\}$.

\noindent
4) We say $\varphi_2(\bar z)$ is complete \when \, for any $\mathbf
s$-structure $M$, if $\bar a_1,\bar a_2 \in \varphi_2(M)$ then $\bar
a_1,\bar a_2$ realizes the same quantifier free type in $M$.

\noindent
5) We say $\bar\varphi$ is complete \when \, $\varphi_2(\bar z)$ and each
$\varphi_{0,i}(\bar x_i,\bar z)$ is (not contradictory and is)  
complete (see (4)) and 
$\varphi_{0,i}(\bar x_i,\bar z) \vdash \varphi_2(\bar z)$.  If not
said otherwise, we assume $\bar\varphi$ is complete.
\end{definition}

\begin{observation}
\label{a7}
1) In Definition \ref{a6}(2), $H_{\bar\varphi,M,\bar c}$ is indeed a graph
(possibly empty) and is finite \when \, 
$M$ is finite $\tau_{\mathbf s}$-structure.

\noindent
2) For each $\bar\varphi$ as in \ref{a6}(1), 
for each $i < \mathbf
i(\bar\varphi)$ one of the following holds:
\mn
\begin{enumerate}
\item[$(\alpha)$]  for some $k,\varphi_2(\bar z) \vdash (\exists!^k
  \bar x_i)(\varphi_i(\bar x_i,\bar z))$
\sn
\item[$(\beta)$]  for every $k$ for some $\mathbf s$-structures $M$, in
  $M$ we have $\varphi_2(\bar z) \vdash (\exists^{\ge k} \bar x)
\varphi_{0,i}(\bar x_i,\bar z)$.
\end{enumerate}
\end{observation}

\begin{PROOF}{\ref{a7}}
Read Definition \ref{a6}(1).
\end{PROOF}

\begin{observation}
\label{a8}
1) Let $\mathbf s$ be an $I$-kind and $\bar\varphi$ is a complete 
$\mathbf s$-scheme. 

The following are equivalent:
\mn
\begin{enumerate}
\item[$(a)$]  for every $\bar p \in \mathbf P^2_{\mathbf s}$ and random
  enough $\cM = \cM_{\mathbf s,n}$ we have $\varphi_2(\cM) \ne \emptyset$
\sn
\item[$(b)$]  for some $\bar p \in \mathbf P^2_{\mathbf s} \cup 
\mathbf P^0_s$ we have $0 < \lim \sup_n 
\Prob(\varphi_2(\cM_{\mathbf s,p,n}) \ne \emptyset)$.
\end{enumerate}
\mn
2) For any sentence $\psi \in \bbL(\tau_{\mathbf s})$, 
similarly replacing $\varphi_2(\cM) \ne \emptyset$ by ``for some $\bar c,
H_{\bar\varphi,\cM,\bar c} \models \psi$".
\end{observation}

\begin{PROOF}{\ref{a8}}
Easy.
\end{PROOF}

\begin{definition}  
\label{a9}
1) We call an $\mathbf s$-scheme $\bar\varphi$ trivial \when \, for each
 $i < \mathbf i(\varphi)$ we have $\ell g(\bar x_i)=0$.

\noindent
2) We call an $\mathbf s$-scheme $\bar\varphi$ degenerated \when \, the
conditions of \ref{a8} fail; as long as $\bar\varphi$ is
complete this does not occur 
as $ \mathbf{i} ( \bar{ \varphi } ) \not= 0$, \medskip Def \ref{a6}(1)$ ab _3 $.

\noindent
3) We say the $\mathbf s$-scheme $\bar\varphi$ is 1-weak \when \, at least one
of the following holds:
\mn
\begin{enumerate}
\item[$(a)$]  $\mathbf s$ is degenerated \oor \, $\mathbf s$ is trivial, 
i.e. $\ell g(\bar x_i)=0$ for every 
$i < \mathbf i(\varphi)$ \oor \, 
\sn
\item[$(b)$]  for some truth value $\mathbf t$ and $i_1,i_2 < \mathbf
  i(\varphi)$ satisfying $\ell g(\bar x_{i_1}),\ell g(\bar x_{i_2})
  \ge 1$ and $v_1 \subsetneqq \ell g(\bar x_{i_1}),v_2 \subsetneqq
  \ell g(\bar x_{i_2})$ we have
\sn
\begin{enumerate}
\item[$\bullet$]  for some (equivalently any) $\bar p \in \mathbf P^2_{\mathbf s}$,
for random enough $\cM = \cM_{\mathbf s,\bar p,n}$, for some $\bar c \in 
\varphi_2(\cM)$ and $\bar a^*_\ell \in \varphi_{1,i_\ell}(\cM,\bar c)$
for $\ell=1,2$ we have
\sn
\item[$\bullet$]  if $\bar a_\ell \in \varphi_{1,i_\ell}(\cM,\bar c)$
  and $\bar a_\ell \rest v_\ell = \bar a^*_\ell \rest v_\ell$ for
$\ell=1,2$ and $\rang(\bar a_1) \cap \rang(a_2) \subseteq \rang(\bar
  a^*_1 \rest v_1) \cap \rang(\bar a^*_2 \rest v_2)$
\then \, $\cM \models \varphi_{i_1,i_2}[\bar a_1,
\bar a_2,\bar c]^{\iif(\mathbf t)}$.
\end{enumerate}
\end{enumerate}
\mn
4) We say the $\mathbf s$-scheme $\bar\varphi$ is 2-weak \when \, at least
one of the following holds:
\mn
\begin{enumerate}
\item[$(a)$]   it is degenerated \oor \, 
trivial, i.e. as in (a) of part (3)
\sn
\item[$(b)$]  for some $i < \mathbf i(\varphi),\ell g(\bar x_i) \ge 2$
\sn
\item[$(c)$]  for some $i_1,i_2 < i(\varphi)$ with 
$\ell g(\bar x_{i_2}) = 1 =
\ell g(\bar x_{i_2})$ and $\bar p \in \mathbf P^2_{\mathbf s}$ and
random enough $\cM = \cM_{s,\bar p,n}$ and $\bar c \in
\varphi_2(\cM)$ there is $\mathbf t \in \{0,1\}$ such that for every
$a_1 \in \varphi_{1,i_1}(\cM),a_2 \in \varphi_{1,i_2}(\cM)$ we have
$a_1 \ne a_2 \Rightarrow H_{\bar\varphi,\cM,\bar c} 
\models ``a_1 R a_2$ iff $\mathbf t=1"$.
\end{enumerate}
\mn
5) We say the $\mathbf s$-scheme is 3-weak \when \, it is 1-weak or 2-weak.
\end{definition}

\begin{claim}  
\label{a15}
1) For any $k$, if $\cM = \cM_{\mathbf s,\bar p,n}$ 
is random enough for $k$ and $\bar c \in {}^{k \ge} M$, and 
there is an interpretation using as parameter
$\bar c \in {}^{k \ge}\cM$ of a graph $H$ in $\cM$
using $(\le k)$-tuples (in the widest sense - the elements can be
equivalence classes of suitable definable equivalence 
relations on set of tuples
satisfying a formula) by formulas of length $\le k$ \then \, there is
a complete 
$\mathbf s$-scheme $\bar\varphi$ such that $H = 
H_{\bar\varphi,\cM,\bar c}$ and $k(\bar\varphi) \le k$.

\noindent
1A)  
For any interpretation by first  order formulas  
 with parameter $ \bar{ z } $  
 \begin{enumerate} 
\item[(*)]  there is an $ \mathbf{s} $-scheme
    interpretation equivalent to it, and we can compute it, 
\item[(*)]  moreover we can compute a finite sequence
    $ \langle \bar{ \varphi } _i: i < i_* \rangle $,
    each $ \bar{ \varphi } _i $ is complete, all with the same parameter $ \bar{ c } $.
\end{enumerate} 

\noindent
2) In fact $\bar\varphi$ depends just on the interpretation and the
 quantifier free type of $\bar c$ in $\cM$, not on $\cM$ (and even $n$).

\noindent
3) Given $\mathbf s$ and $k$ there only finitely many 
scheme $\bar\varphi$ as above.  
\end{claim}

\begin{PROOF}{\ref{a15}}  
Obvious.
\end{PROOF}

\begin{definition}
\label{a18}
Let $\mathbf s,\bar\varphi$ be as above,  
    $ \bar{ \varphi } $ is complete, see \ref{a6}(5). 

\noindent 
We say $(\mathbf s,\bar\varphi)$ is reduced \when \,: for every $\bar
p \in \mathbf P^2_{\mathbf s}$ and random
enough $\cM = \cM_{\mathbf s,\bar p,n}$ and $\bar c \in {}^{\ell
  g(\bar z_{\bar\varphi})}\cM$ satisfying $\varphi_2(\bar
z_{\bar\varphi})$, the graph $H = H_{\bar\varphi,\cM,\bar c}$ is not
$H = H_{\bar\varphi',\cM,\bar c'}$ when $(\bar\varphi',\bar c'$
appropriate and) $\Rang(\bar c') \subsetneqq \Rang(\bar c)$; recall $\bar
c$ is without repetitions.

\end{definition}
\bigskip

\subsection {Simple Random Graph}\
\bigskip

Our intention is that the behaviour of $\cG_{q,n } $  
expanded by
some formulas in the expanded logic will be like $\cM_{\mathbf s,\bar
  p , n},\bar p \in \mathbf P^2_{\mathbf s}$, but we need a relative 
  as we can iterate.  

\begin{definition}
\label{a21}
For $\iota = 1,2,3$ let $\mathbf U_\iota$ be the set of 
objects $\mathbf u$ consisting of the following 
(we may add subscript $\mathbf u$):
\mn
\begin{enumerate}
\item[$(a)$]  $\bar{\mathbf s} = \langle \mathbf s_\ell:\ell \le
  \ell(\mathbf u)\rangle$
\sn
\item[$(b)$]  $\mathbf s_\ell$ is a kind sequence
\sn
\item[$(c)$]  $\mathbf s_0 = \mathbf s_{\gr}$, the graph kind sequence,
  see \ref{a2}(1A)
\sn
\item[$(d)$]  $\mathbf s_\ell \subseteq \mathbf s_{\ell +1}$,
  i.e. $I_{\mathbf s_\ell} \subseteq I_{\mathbf s_{\ell +1}}$ and $t \in
  I_{\mathbf s_\ell} \Rightarrow (k_{\mathbf s_\ell,t},K_{\mathbf s_\ell,t})
  = (k_{\mathbf s_{\ell +1},t},K_{\mathbf s_{\ell +1},t})$
\sn
\item[$(e)$]  notation: so we may write $(k_{\mathbf u,t},K_{\mathbf u,t})$
  for $t \in I_{\mathbf s_{\ell(\mathbf u)}}$ and $I_{\mathbf q} = 
I_{\mathbf s_{\ell(\mathbf u)}}$
\sn
\item[$(f)$]  for $t \in I_{\mathbf s_{\ell +1}} \backslash 
I_{\mathbf s_\ell}$ we have: $\bar\varphi_t$
is 
a complete reduced
$\mathbf s_\ell$-scheme, not $\iota$-weak such that
$K_t = K_{\bar\varphi_t}$, see Definition \ref{a45}(2) let $\mathbf i_t
= \mathbf i(t) = \mathbf i(\bar\varphi_t)$ and similarly
$\varphi_{t,2},\varphi_{t,0,i},\varphi_{t,1,i,j}$ but let
$\varphi_t(\bar z_t) = \varphi_{t,2}(\bar z_t)$.
 In the case $\iota = 2, 3 $ 
 if $\bar y_{t,i} \ne \langle \rangle$ then 
$\bar y_{t,i}$ is a singleton so we shall write $\varphi_{t,0,i}
(y,\bar z_{t,i})$
\sn
\item[$(g)$]  $q = q_{\mathbf u} \in (0,1)_{\bbR}$.
\end{enumerate}
\end{definition}

\begin{definition}
\label{a24}
For $\mathbf u \in \mathbf U_\iota$ we define a random $\cM_{\mathbf u,n}$,
i.e. a 0-1 context, as follows.

For a given $n,\cM_{\mathbf u,n}$ is gotten by drawing $\cM_{\mathbf u,n,\ell}
\in \mathbf M_{\mathbf s_{\mathbf u,\ell},n}$ by induction on $\ell \le
\ell(\mathbf u)$ and in the end $\cM_{\mathbf u,n} = {\mathscr M }  
_{\mathbf
  u,n,\ell(\mathbf u)}$.

Now
\mn
\begin{enumerate}
\item[$(a)$]  if $\ell=0,\cM_{\mathbf u,n,\ell}$ is $\mathbf 
{\mathscr G}  
_{q(\mathbf
    u),n}$, i.e.the random graph on $n$ with edge probability $q$
\sn
\item[$(b)$]  if $\ell < \ell(\mathbf u)$ and $\cM_{\mathbf u,n,\ell}$ has
  been drawn and $t \in I_{\mathbf s_{ell +1}} \backslash I_{\mathbf
    s_\ell}$, we draws $R_t(\cM_{\mathbf s_{\ell +1}})$ as follows:
\sn
\begin{enumerate}
\item[$(\alpha)$]  if $\bar c \in \varphi_t(M)$ we draw the truth
  value of $\bar c \in R_t(\cM_{\mathbf s_{\ell +1},n})$ with
  probability $\mathbf h(\sum\limits_{i<\mathbf i(t)} 
  \EXP |\varphi_{a,i}(\cM_{\mathbf
    s_{\ell +1},n},\bar c)|/|K_{t,i}|)$
    recalling $ \EXP $ is the expected value
\sn
\item[$(\beta)$]  if $\bar c$ is a sequence of length $k_t$ but
$\notin \varphi_t(M)$ then $\bar c \notin R_t(\cM_{\mathbf s_{\ell +1},t})$.
\end{enumerate}
\end{enumerate}
\end{definition}

\begin{claim}
\label{a27}
For $\mathbf u \in \mathbf U_\iota,\cM_{\mathbf u,n}$ is like $\cM_{\mathbf s_{\mathbf
    q},\bar p}$ for any $\bar p \in \mathbf P^2_{\mathbf s_{\mathbf q}}$
(and $M_{\mathbf u,n,\ell}$ like $\cM_{\mathbf s_{\mathbf q,\ell},\bar p})$,
in particular, satisfying the zero one law:
\mn
\begin{enumerate}
\item[$(*)$]  for any $k_1$ for some $k_2$, for any random enough
  $\cM_{\mathbf u,n}$ we have:
\sn
\begin{enumerate}
\item[$\bullet$]  if $\varphi(\bar x),\psi(\bar y,\bar z)$ are
  complete $\bbL(\tau_{\mathbf s_{\mathbf u}})$-formulas such that
$\psi(\bar y,\bar z) \vdash \varphi(\bar z)$ (so they respect the
$K_{\mathbf u,t}$'s!, see Definition \ref{a9}(6)) 
and $\ell g(\bar y)+ \ell g(\bar x) \le
k_1$ and $\bar c \in \varphi(\cM_{\mathbf u,n})$ and $k_{t,i} \ge 1$
\then\, the number of members of $\psi_{t,i}(M_{\mathbf u,n},\bar c)$ is
similar to $\binom {n^{\ell g(\bar y)}}{k_t} \cdot \frac{k_t}{(K_t)}$;
fully
\sn
\item[$\bullet$]  at most\footnote{we could have allowed, e.g. when
    $k_t=1$ to be near to 1 though not too closely, but if we shal use
    a quantifier $\mathbf Q$ such that $\ll \frac 12$ of the structures
    are in it}
$\binom {n^{\ell g(\bar y_{t,i})}}{k_{t,i}} \cdot
\frac{k_{t,i}!}{|K_{t,i}|} \cdot (1 + \frac{1}{k_2})$
\sn
\item[$\bullet$]  at least
$\binom {n^{k_{t,i}}}{k_t} \cdot \frac{k_{t,i}!}{|K_{t,i}|} \cdot
\mathbf h(n) - k_2$
\sn
\item[$\bullet$]  if $\iota=2$, then $k_{t,i}=1$, so this is simpler.
\end{enumerate}
\end{enumerate}
\end{claim}

\begin{remark}  \label{a25}
What is the reason for our choice in Clause $(b)( \alpha ) $ 
of Def \ref{a24}? There are some demands pulling in different directions.
\begin{enumerate} 
\item[(a)]  This probability should  be not too small (considering 
it belongs to $ (0, 1/2)$) such that the argument 
\lqq a $  \Sigma _1 $ formulas $ (\exists  y) \varphi ( y, \bar{ a } ) $ 
    hold when not excluded" as in $ {\mathscr M } _{\mathbf{s} , \bar{ p }, n }$ 
\item[(b)] but always is not so small such that
  $ \Prob((\exists  \bar{ y } )\varphi ( \bar{ y }  \varphi )$ 
converge to zero or to one
\item[(c)] The $ {\mathscr M } _{\mathbf{u} , n}$ are intended
to imitate what we get by starting with ${\mathscr G} _{p,n} $
    and expanding it by relations definable by formulas 
    $ \varphi(  \bar{ x } )$ from our logic, 
    so we are applying our quantifier to a definable (with parameters)
    graph. So such a graph even  almost surely will not have exactly
    $ n $ nodes. In the non-degenerated case the number will be 
    of the order of magnitude 
    \begin{enumerate} 
\item[(*)] $ C n^k $ for some positive real $ C $ and $ k \ge 1$
in the $ 1$-low case   
\item[(*)] $ C n $ for some positive real $ C $
 in the $ 2/3$-low case,   
    \end{enumerate} 
\end{enumerate} 
\end{remark} 

\begin{PROOF}{\ref{a27}}
Should be clear.
\end{PROOF}
\bigskip

\subsection {Low/High Graphs}\
\bigskip

\noindent
An $\mathbf s$'s scheme $\bar\varphi$ may be such that, e.g. the 
bi-partite graph with the $i$-th
kind and the $j$-th kind is 
    in the low case, see 
Definition \ref{a9}(4);  so 
we try
to single out those $\bar\varphi$'s.  Those cases are ``undesirable"
for us and we shall try to discard them.

\begin{definition}
\label{a34}
1) We say a finite graph $H$ is $\mathbf h-1$-low (recall $\mathbf h$ is
from \ref{a0} so can be omitted)  \when \, there are no disjoint
$A,B \subseteq H$ and $\iota < 2$ such that (letting $n=|H|$)
\mn
\begin{enumerate}
\item[$(a)$]  $|A|,|B| \ge |H|^{\mathbf h(n)}$
\sn
\item[$(b)$]  if $a \in A$ and $b \in B$ \then \, $(a,b)$ is an edge of
$H$ \Iff \, $\iota=1$.
\end{enumerate} 
\mn
2) We say that a finite graph $H$ is $\mathbf h-2$-low \when
\,\footnote{The specific choice of $m$ is not important, but they have
  to be $\le n^{1/k}$ and $>k$ for any $k$, for large enough $n$.
  Similarly $|\Rang(\mathbf c_\ell)|$ compared to $m$.}
 letting $n = |H|,m = \lfloor \mathbf g(n) \rfloor = \lfloor n^{\mathbf
   h(n)} \rfloor$, there are no $\bar a,\bar b,M,\mathbf c$ such that:
\mn
\begin{enumerate}
\item[$(a)$]  $\bar a = \langle a_\ell:\ell < m\rangle$
\sn
\item[$(b)$]  $\bar b = \langle b_{\ell,k}:\ell < k \le m\rangle$
\sn
\item[$(c)$]  $\bar a \char 94 \bar b$ is a sequence of nodes of $H$
  with no repetitions 
\sn
\item[$(d)$] each  
    $\mathbf c_0,\mathbf c_1$ is a function from
$\{(\ell,j):\ell,j \le m\}$ to $\{0,1,\dotsc,\lfloor
\mathbf g(n) \rfloor\}$
\sn
\item[$(e)$]  $\mathbf c_2$ is a function from $\{(\ell,k):\ell,k \le m\}$ 
into $\{0,1,\dotsc,\lfloor \mathbf g(n) \rfloor\}$
\sn
\item[$(f)$]  if $\ell' < k' \le m$ and $j'<m$ and $\ell'' < k'' \le
m,j'' \le n$ and $\mathbf c_0(\ell',j') = \mathbf c_0(\ell'',j'')$
and $\mathbf c_1(k',j') = \mathbf c_1(k'',j'')$
 and $\mathbf c_2(k',\ell) = \mathbf c_2(b'',\ell'')$ 
\then \, $(b_{\ell',k'},a_{j'})$ is an edge of $H$ iff
$(b_{\ell'',k''},a_{j''})$ is an edge of $H$.
\end{enumerate}
\mn
3) We say the finite graph $H$ is $\mathbf h-3$-low when it is $\mathbf
h-1$-low or $\mathbf h-2$-low.

\noindent
4) In parts 
(1) and (2), $\mathbf h-\iota$-high means the negation of $\mathbf
h-\iota$-low. 
\end{definition}

\begin{claim}
\label{a37}
Assume $\mathbf s$ is an $I$-kind, (see Definition \ref{a2}) and
$\bar\varphi$ is a complete   
$\mathbf s$-scheme (see Definition\ref{a9}(2)).
\ref{a6}, \ref{a9}(2))
\mn
\begin{enumerate}
\item[$(A)$]  the following are equivalent:
\sn
\begin{enumerate}
\item[$(\alpha)$]  $\bar\varphi$ is trivial
\sn
\item[$(\beta)$]  if $\bar p \in \mathbf P^2_{\mathbf s}$ \then \,
for random enough $\cM = \cM_{\mathbf s,n,\bar p}$ and
  $\bar c \in \varphi_2(\cM)$ the graph $H_{\bar\varphi,\cM,\bar c}$
  has $\le \mathbf i(\bar\varphi)(k(\bar\varphi)!)$ nodes
\sn
\item[$(\gamma)$]  if $\varepsilon > 0$ and $\bar p \in \mathbf
  P^2_{\mathbf s}$ \then \, $0 < \lim\sup_n \Prob$(letting $\cM = 
\cM_{\mathbf s,\bar p,n}$, for some $\bar c \in \varphi_2(\cM)$ the 
graph $H_{\bar\varphi,\cM,\bar c}$ has $\le n^{1-\varepsilon}$ nodes).
\end{enumerate}
\sn
\item[$(B)$]  the following are equivalent for non-trivial $\bar\varphi$:
\sn
\begin{enumerate}
\item[$(\alpha)$]  $\bar\varphi$ is 1-, see Def \ref{a9},
\sn
\item[$(\beta)$]   if $\bar p \in \mathbf P^2_{\mathbf s}$ \then \,
for every random enough $\cM = \cM_{\mathbf s,n,\bar p}$
  and for every $\bar c \in \varphi_2(\cM)$ the graph 
$H_{\bar\varphi,\cM,\bar c}$ is $\mathbf h-1$-low
\sn
\item[$(\gamma)$]  if $\varepsilon > 0$ and $\bar p \in \mathbf
  P^2_{\mathbf s}$ \then \, $0 < \lim\sup_n \Prob$(letting $\cM = 
\cM_{\mathbf s,\bar p,n}$, for some $\bar c \in \varphi_2(\cM)$ the 
graph $H_{\bar\varphi,\cM,\bar c}$ is 1-low)
\end{enumerate} 
\sn
\item[$(C)$]  Like (B), replacing 1-weak, $\mathbf h-1$-low by 2-weak,
  $\mathbf h-2$-low respectively
\sn
\item[$(D)$]  Like (B), replacing 1-weak, $\mathbf h-1$-low by 3-weak,
 $\mathbf h-3$-low respectively.
\end{enumerate} 
\end{claim}

\begin{PROOF}{\ref{a37}}
\medskip

\noindent
\underline{Clause (A)}:

Trivially $(A)(\alpha) \Rightarrow (A)(\beta)$ and $(A)(\beta)
\Rightarrow (A)(\gamma)$.

So it suffices to assume $\bar\varphi$ is non-trivial, $\bar p \in
\mathbf P^2_{\mathbf s}$ and let $\varepsilon > 0$ be small enough and
prove that for every random enough $\cM = \cM_{\mathbf s,\bar p,n}$ and
$\bar c \in \varphi_2(\cM)$ the graph $H_{\bar\varphi,\cM,\bar c}$ has
$\ge \varepsilon n$ nodes.

Let $i < \mathbf i(\bar\varphi)$ be such that $k_i = \ell g(\bar x_i)
>0$, so for $n$ large enough and $\bar c \subseteq [n]$ of length
$\ell g(\bar z)$ let $S_{n,\bar c} = \{\bar a:\bar a$ is a
sequence of length $\ell g(\bar x_i)$ with no repetition of members of
$[n]$ not from $\bar c\}$.  For every $\bar a \in S_{n,i}$, the real
$\Prob(\cM_{\mathbf s,\bar p,n} \models$ ``if $\varphi_2(\bar c)$
then $\varphi_{1,i}(\bar a,\bar c)")$ is the same for every $\bar a \in
S_{n,\bar c}$ and is of the form $r(1) \mathbf g(n)^m$ for some $r(1) \in
(0,1)_{\bbR},m \in \bbN \backslash \{0\}$ not depending on $n$.
Fixing $\bar c$ under the assumption $\cM_{\mathbf s,\bar p,n} \models
\varphi_2[\bar c]$, considering a maximal set of pairwise disjoint $i
\in S_{n,\bar c}$, the events $\cM_{\mathbf x,\bar p,n} \models
\varphi_{1,i}[\bar a,\bar c]$ are independent, such that almost surely the
number $|\{\bar a \in S_{n,\bar c}:\cM_{\mathbf s,\bar p,n} \models
\varphi_{1,i}[\bar a,\bar c)\}$ is $\ge n/(r(1) 
\mathbf g(n)^m(1-\varp))$.  Similarly almost surely the number of $\bar c$
  such that $\cM \models \varphi_2[\bar c]$ is large.
\medskip

\noindent
\underline{Clause (B)}:

\underline{First why $(B)(\alpha) \Rightarrow (B)(\beta)$}?

Note $\bar\varphi$ is non-trivial; 
$(\mathbf s,\bar\varphi)$ cannot satisfy clause (a) of Definition
\ref{a9} because in the present claim we are assuming $\bar\varphi$ is
non-degenerated.   So assume clause (b) of \ref{a9}(3) holds as
exemplified by $i_1,i_2,\mathbf i(\bar\varphi),v_1,v_2$ and truth value
$\mathbf t$, i.e. $\ell g(\bar x_i),\ell g(\bar x_j) > 0$, etc.  So assume
$n$ is large enough and $\cM = M,\bar c \subseteq [n]$ has length
$\ell g(\bar z_{\bar\varphi})$.

Let $A_\ell = \{\bar a:\bar a \subseteq [n]$ is of length $\ell g(\bar
x_{i_\ell})$ for $\ell=1,2$ with no repetition and is disjoint to
$\bar c\}$.  Choose for $\ell=1,2$ disjoint
$\bar a^*_\ell \in A_\ell$.  So the event $\cE_{\bar c} =
``(\bar c \char 94 \bar a^*_1 \char 94 \bar a^*_n)$ is as in
\ref{a9}(3)" has probability $\ge r_1(\mathbf g(n)^{k(1)}) $ for some $r_1
\in (0,1)_{\bbR},k \in \bbN \backslash \{0\}$ not depending on $n$
(and $\bar c$).  Fixing $(\bar c,a^*_1,a^*_2)$ let $C_\ell \subseteq
(n^1 \backslash \rang(\bar c \char 94 \bar a^*_1 \char 94 \bar
a^*_2),|C_\ell| \ge (n-|\ell g|(\bar z \char 94 \bar x_{i_1} \char 94
\bar x_{i_1}) - \frac 12$ for $\ell=1,2$ and $C_1 \cap C_2 =
\emptyset$.  Let $A'_\ell = \{\bar a \in A_\ell:\Rang(\bar a)
\subseteq C_\ell\}$.

Easily for some $r(2),r(3) \in (0,1)_{\bbR}$ not depending on $n,t$
 the probability of the event $\cE_2 = \cE^r_{\bar c,\bar a^*,\bar a^*_1}$
 is $\ge 1 - 2^{-r(2)n}$ where
\mn
\begin{enumerate}
\item[$(*)$]  $\cE_2$ means: if $\cM \models \varphi_2[\bar c] \wedge
  \varphi_{1,i_1}[\bar a^*_1] \wedge \varphi_{1,i_2}[\bar a^*_2]$
  \then \, $|\{\bar a_\ell \in A'_\ell:\cM \models \varphi[\bar
  a_{\ell,m}]\}| \ge n^{|u(\ell)|}r(3)$ for $\ell=1,2$.
\end{enumerate}
\mn
If $\cE_2$ occurs, clearly $\mathbf t$ and $A^*_{\cM,\ell} = 
\{\bar a/E_{K_{\mathbf s,i,\ell}}:\bar a \in A'_\ell$ and $\cM \models
\varphi_{0,i_\ell}[a_{n,m},\bar c]\}$ for $\ell=1,2$ exemplifies
$H_{\bar\varphi,\cM,\bar c}$ is low.  As the number of $\bar c,\bar
a^*_1,\bar a^*_2$ is polynomial we can finish.
\medskip

\noindent
\underline{Second, why $(B)(\beta) \Rightarrow (B)(\gamma)$}:

Read the clauses and Definition of \ref{a34}.
\medskip

\noindent
\underline{Third, $\neg(B)(\alpha) \Rightarrow \neg(B)(\gamma)$}:
This suffices

Why this holds?  Let $\cM = \cM_{\mathbf s,\bar p,n}$ be random enough,
$\bar c_2 \in \varphi_2(\cM)$ and $A_1,A_2 \subseteq H =
H_{\bar\varphi,\cM,\bar c}$ witness $H$ is low, so $|A_\ell| \ge
n^{\mathbf h(n)}$.  So $n^*_1 = \min\{|A^*_1|;|A^*_2|\} \ge m^{\mathbf h(n)}$.

Clearly for each $\ell \in \{1,2\}$ for some $i(\ell) < \mathbf
i(\bar\varphi)$ we have

\[
|\bar a/K_{\mathbf s,i,\ell} \in A_\ell:\bar a \in
\varphi_{1,i}(\cM,\bar c)\}| \ge |A_\ell|/\mathbf i(\bar\varphi) = n^*_2
\ge n^{\mathbf h(n)}/\mathbf i(\bar\varphi).
\]

\mn
So for some $r \in (0,1)_{\bbR}$ not depending on $n$
 for $\ell=1,2$ we can find $\langle \bar a_{\ell,m}:m < n^*_3 =
 \langle (n^*_2)^r \rangle$ and partition $v_\ell,u_\ell$ of $\ell
 g(\bar x_{i(\ell)})$ such that:
\mn
\begin{enumerate}
\item[$(*)$]
\begin{enumerate}
\item[(a)]  $\bar a_{\ell,m} \rest v_c = a^*_\ell$
\sn
\item[(b)]  $\Rang(\bar a_{\ell_1,m_1} \rest u_\ell) \cap
\Rang(\bar a_{\ell_2,m_2}) = \emptyset$ when $m_1,m_2 <
n^*_3(\bar\varphi) \wedge \ell_1,\ell_2 \in \{1,2\} \wedge (\ell_1,m_1)
\ne (\ell_2,m_2)$
\sn
\item[(c)]  $\Rang(\bar a_{\ell,m} \rest u_\ell),
\Rang(\bar a_{2,m(r)} \rest u_2),\bar a^*_1 \char 94 \bar a^*_2$ are
pairwise disjoint for $\ell \in \{1,2\},m < n^*_3$.
\end{enumerate}
\end{enumerate}
\mn
We draw $\cM \rest (\bar c \char 94 \bar a_{\ell,m})$ for every $\ell
\in \{1,2\}$ and $m < n^*_3$ we get $\cM'$.  So ignoring events of very
low probability ($\le ({\frac{1}{2})^{rn}}$ for fix $r \in (0,1)_{\bbR}$)
\mn
\begin{enumerate}
\item[$(*)$]  $w_\ell := \{m < n^*_3:(M' \rest \bar c \char
  94 \bar a_{\ell,m}) \models \varphi_{1,i(\ell)}[\bar a_{\ell,m},\bar
  c]\}$ has $\ge n^*_4 := \sqrt n^*_3$ members.
\end{enumerate}
\mn
So $n^*_4 \ge n^\varepsilon$ for $\varepsilon$ small enough but let $Y_\ell =
\{\bar a_{\ell,m}/K_{t,i(\ell)}:m \in w_\ell\}$; it is a set of $\ge
n^*$ nodes of $H_{\bar\varphi,\cM,\bar c}$.

Now
\mn
\begin{enumerate}
\item[$(a)$]  $m(1),m(2) < n^*_4 \Rightarrow 
\Prob(\cM \models \varphi_{1,i(1),i(2)}(\bar a_{1,m(1)},
\bar a_{2,m(2)},\bar c)) = r/\mathbf g(n)^k$ for some $r 
\in \bbR_{>0},k \in \bbN \backslash \{0\}$ not depending on $n$
\sn
\item[$(b)$]  if $i(1),i(2)$ are not as required in \ref{a9}(3)(h) and
  $\mathbf t =0,1$ then with negligible probability we have for some
  $u_1 \subseteq w_2,u_2 \subseteq w_2$ with $\lfloor \mathbf g(|H|)
  \rfloor$ elements each we have $m(1) \in u_1 \wedge m(2) \in u_2
 \Rightarrow \cM \models \varphi_{1,i(1),i(2)}
(\bar a_{1,m(1)},\bar a_{2,m(2)},\bar c)^{\iif(\mathbf t)}$.
\end{enumerate}
\mn
So this could not have occured.
\medskip

\noindent
\underline{Clauses (C),(D)}:

Also straightforward.
\end{PROOF}
\bigskip

\centerline {$*\qquad * \qquad *$}
\bigskip

\begin{definition}
\label{a45}
1) Assume $\bar\varphi^1 = \bar\varphi_1,\bar\varphi^2 = \bar\varphi_2$
are $\mathbf s$-schemes and 
 $\bar\varphi_1,\bar\varphi_2$ are reduced and complete.  We say $(\mathbf
s,\bar\varphi^1),(\mathbf s,\bar\varphi^2)$ are explicitly isomorphic
\when \, some $\pi$ and $\varkappa$ witness it which means:
\mn
\begin{enumerate}
\item[$(a)$]  $\mathbf i(\bar\varphi^1) = \mathbf i(\bar\varphi^2)$ and
$\ell g(\bar z_{\bar\varphi^1}) = \ell g(\bar z_{\bar\varphi^2})$
\sn
\item[$(b)$]  $\pi$ is a permutation of $\{0,\dotsc,\mathbf
  i(\bar\varphi^1)-1\}$ such that $k_{\bar\varphi^1,i} =
k_{\bar\varphi^2,\pi(i)}$ and $K_{\bar\varphi_1,i} =
K_{\bar\varphi_2,i}$ for $i < \bar\varphi_1$
\sn
\item[$(c)$]  $\varkappa$ is a permutation of 
$\ell g(\bar z_{\bar\varphi_1})$
\sn
\item[$(d)$]   for random enough $\cM = \cM_{\mathbf s,\bar p,n}$, if
  $\ell \in \{1,2\},\cM \models \varphi_2^\ell[\bar c_\ell]$ \then\,
  letting $\bar c_{3 -\ell}$ be such that $\bar c_2 = \varkappa(\bar c_1)$
  we have $M \models \varphi^{3-\ell}_2[\bar c_{3-\ell}]$ and
$\varphi_{1,i}(\cM,\bar c_1) = \varphi_{1,\pi(i)}(\cM,\bar c_2)$ and
$\varphi_{1,i,j}(\cM,\bar c_1) = \varphi_{1,\pi(i),\pi(j)},(\cM,\bar c_2)$
\end{enumerate}
\mn
2) For $\mathbf s,\bar\varphi$ as above let $K_{\bar\varphi} = K_{\mathbf
s,\bar\varphi}$ be the group of permutations $K$ of $\ell g(\bar
z_{\bar\varphi})$ such that $\bar\varphi$ is explicitly isomorphic to
itself using our $\varkappa$ in \ref{a45}(1).
\end{definition}

\begin{claim}
\label{a48}
1) For every $\mathbf s$-scheme $\bar\varphi$ we can find $\langle
\bar\varphi^\iota(\bar z_\iota):\iota < \iota(*)\rangle$ such that:
\mn
\begin{enumerate}
\item[$(a)$]  $\bar\varphi^\iota(\bar z^\iota_i)$ is a complete
reduced $\mathbf s$-scheme such that $\bar z^\iota$ is a subsequence of
$\bar z$
\sn
\item[$(b)$]  for every $\mathbf s$-structure $M$ and $\bar c 
\in \varphi_2(M)$ for some
$\iota$ letting $\bar c_\iota = \langle c_j:j \in \dom(\bar z)$ and
$z_j$ appears in $\bar z^\iota\}$ we have $H_{\bar\varphi,M,\bar c} \cong
H_{\bar\varphi^\iota,M,\bar c^\iota}$
\sn
\item[$(c)$]  for every $\mathbf s$-structure $M$
 $\iota < \iota(*)$ and $\bar c^\iota \in
\varphi^\iota_2(M)$ there is $\bar c$ such that $(\bar c,\bar
c^\iota,\bar\varphi,\bar\varphi^i)$ are as in clause (b).
\end{enumerate}
\mn
2) For complete $\bar\varphi$ in the definition of ``trivial",
``degenerated", ``reduced" we can replace ``some $\bar c$" by ``
$\bar c'$".

\noindent
3) In the definition of $\bbL(\mathbf Q_{\mathbf t})(\tau)$, see
Definition \ref{b2}, we can use $(\mathbf Q \dotsc,\bar x_{1,i},\bar
   x'_{1,i},\ldots) \bar\varphi$ for complete reduced
   non-trivial, non-degenerated $\bar\varphi$.
 \end{claim}

\begin{PROOF}{\ref{a48}}
Easy.
\end{PROOF}

\begin{tic}
\label{a51}
Assume $\mathbf s$ is an $I$-kind and 
$\bar\varphi',\bar\varphi''$ are complete reduced $\mathbf s$-schemes as above.

\noindent
1) If $\cM = \cM_{\mathbf s,\bar p,n}$ is random enough
and $\cM \models \varphi'_2[\bar c'] \wedge \varphi''_2[\bar c'']$ so
$H' = H_{\bar\varphi',\cM,\bar c'},H'' = H_{\bar\varphi'',\cM,\bar c''}$
are well defined \then \, $H' \cong H''$ \Iff \, $\Rang(\bar
c'),\Rang(\bar c'')$ and moreover
$(\mathbf s,\bar\varphi'),(\mathbf x,\bar\varphi'')$ are explicitly
isomorphic, as witness by $(\pi,\varkappa)$ such that $\pi$ maps
$\bar c'$ to $\bar c''$, see Definition \ref{a45}.

\noindent
2) Being explicitly isomorphic $\mathbf s$-schemes is an equivalence relation.
\end{tic}

\begin{PROOF}{\ref{a51}}
Straightforward.
\end{PROOF}
\newpage

\section {The random quantifier} \label{2}

\begin{hypothesis}
\label{b0}
Let $\iota \in \{1,3\}$ but $\iota=3$ is simpler and 
large part is
O.K. also for $\iota=2$.
\end{hypothesis}

\begin{definition}
\label{b2}
1) We say $\mathbf Q = \mathbf Q_{\mathbf K}$ is a $\mathbf
   h-\iota$-high-graph quantifier \when \,:
\mn
\begin{enumerate}
\item[$(a)$]  $\mathbf Q$ is a quantifier on finite graphs, i.e. it is a
  class of finite graphs closed under isomorphisms
\sn
\item[$(b)$]  if $H$ is a finite graph and is $\mathbf h-\iota$-low
\then \, $H \notin \mathbf Q$.
\end{enumerate}
\mn
2) We define a probability space on the set of high-graph quantifiers
as follows:  let $\bar H^* = \langle H^*_m:m \in \bbN \rangle$ be a
sequence of pairwise non-isomorphic finite graphs such that each finite
graph is isomorphic to exactly one of them.

For $\iota \in \{1,2,3\}$, we let: 
\mn
\begin{enumerate}
\item[$(a)$]  $\mathbf T = \mathbf T_\iota = \{\bar{\mathbf t}:\bar{\mathbf t} 
= \langle \mathbf t_m:m \in \bbN\rangle,\mathbf t_m$ a truth value, 
$\mathbf t_m = 0$ if $H^*_m$ is $\mathbf h-\iota$-low$\}$
\sn
\item[$(b)$]  we draw the $\mathbf t_m$'s independently, $\mathbf t_m = 0$
if $H^*_m$ is $\mathbf i-\iota$-low and $\mathbf t_m=1$ has probability
$1/\mathbf g(|H^*_m|)$ when $H^*_m$ is not $\mathbf h-\iota$-low
\sn
\item[$(c)$]  Let $\mu_{\mathbf T_\iota}$ be the derived distribution.
\end{enumerate}
\mn
2A) So the probability space is $(\bbB,\mu_{\mathbf T}),\bbB$
is the family of Borel subsets of ${}^{\bbN}2,\mu_{\mathbf T}$ the measure.

\noindent
3) For $\bar{\mathbf t} \in \mathbf T$ let $\mathbf Q^\iota_{\bar{\mathbf t}}$ be
the quantifier $\mathbf Q_{\mathbf K_{\bar{\mathbf t}}},
\mathbf K_{\bar{\mathbf t}} = \{H:H$ a finite graph isomorphic to some
$H^*_m$ such that $\mathbf t_n=1\}$.

\noindent
4) We say $H$ is $\mathbf h-\iota$-high where $H$ is a finite graph
which is not $\mathbf h-\iota$-low.
\end{definition}

\begin{claim}
\label{b6}
For every random enough $\bar{\mathbf t} \in \mathbf T$ 
the following holds.  

\noindent
1) $\mathbf Q_{\mathbf t}$ is a Lindstr\"om quantifier.

\noindent
2) For random enough graph $\cG_{n,p},\mathbf Q_{\bar{\mathbf t}}$ define non-trivial
   quantifier, defining (with parameters) non-first order definable sets.

\noindent
3) More specifically the formula $\psi(x) =$ (the graph restricted to
   $\{y:yRx\}$ belongs to $\mathbf K_{\bar{\mathbf t}}$)
   define in every random enough $\cG_{n,p}$, a set which is not first
   order definable by a formula of length $k$.
\end{claim}

\begin{PROOF}{\ref{b6}}
Straightforward.
\end{PROOF}

\noindent
So
\begin{definition}
\label{b8}  
1)   
The set of formulas $\varphi(\bar x)$ of $\bbL(\mathbf Q_{\mathbf
   t})(\tau_{\mathbf s})$ for a kind sequence $\mathbf s$ is the closure
   of the set of atomic formulas of $\bbL(\tau_{\mathbf s})$ by negation
   $(\psi(\bar x) = \neg \varphi(\bar x))$, conjunction $(\psi(\bar
   x))=\varphi_1(\bar x) \wedge \varphi_2(\bar x))$, existential
   quantification $(\psi(\bar x) = (\exists y)\varphi(\bar x,y))$ and
   applying $\mathbf Q_{\mathbf t},\psi(\bar z) = (\mathbf Q_{\mathbf
   t},\dotsc,\bar x_{0,i},\bar x'_{0,i},\ldots)_{i < \mathbf
   i(\bar\varphi)} \bar\varphi$ where $\bar\varphi$ is an $\mathbf
   s$-scheme of formulas which are already in 
$\bbL(\mathbf Q_{\mathbf t})(\tau_{\mathbf s})$, so as defined in
   \ref{a6}(1) except that now the $\varphi_{\iota,i}$ are not
   necessarily quantifier free formulas from $\bbL(\tau_{\mathbf s})$.

\noindent
2) Satisfaction, i.e. for a (finite)  
$\mathbf s$-structure $M$, formula
   $\varphi(\bar x)$ and sequence $\bar a$ of elements of $M$ of
   length $\ell g(\bar x)$, we define the truth value of $M \models
   \varphi[\bar a]$ by induction on $\varphi$, the new case is when:
\mn
\begin{enumerate}
\item[$\bullet$]  $\varphi(\bar z_{\bar\varphi}) = (\mathbf Q_{\bar
t},\dotsc,\bar x_{0,i},\bar x'_{0,i},\ldots)_{i < \mathbf i(\bar\varphi)}
\bar\varphi$.
\end{enumerate}
\mn
Now $M \models \varphi[\bar c]$ iff $\bar c \in \varphi_2(M)$ and
$H_{\bar\varphi,M,\bar c}$ is isomorphic to some graph from
$\{H^*_m:\mathbf t_m=1\}$. 

\noindent
3) The syntax of $\bbL(\mathbf Q_{\bar{\mathbf t}})$ does not depend on
   $\bar{\mathbf t}$ so may write $\bbL(\mathbf Q)$ that is 
$\bbL(\mathbf Q)(\tau)$ is the relevant set of formulas, but
   the satisfaction depends so we shall write $M \models_{\bar{\mathbf
   t}} \varphi[\bar a]$ for $\bar a$ a sequence from $M$ and formula
   $\varphi(\bar x) \in \bbL(\mathbf Q)$; of course, such that $\ell g(\bar a) =
   \ell g(\bar x)$.
\end{definition}

\begin{theorem}
\label{b13}
1) For any $p \in (0,1)_{\bbR}$ for all but a null set of 
$\bar{\mathbf t} \in \mathbf T$, the random
   graph $\cG_{n,p}$ satisfies the 0-1 law for the logic 
$\bbL(\mathbf Q^\iota_{\bar{\mathbf t}})$, i.e. we may allow to apply $\mathbf
   Q_{\bar{\mathbf t}}$ to definitions as in Definition \ref{a6}, see
   Claim \ref{a15}.


\noindent  
2)The limit theory $T_*$ is decidable modulo an oracle for the random
${\mathbf K}_ \mathbf{t} $.
\end{theorem}

\begin{remark}
\label{b17}
1) Of course, we can replace the class of graphs by the class of $\mathbf
s$-structures, $\mathbf s$ any kind sequence.

2)  Does the limit theory  depend on $\bar{\mathbf t}$? The problem 
is for when we apply the quantifier . to graphs of fixed size, so
use completer $ \bar{ \varphi } $ with $ k^*_\mathbf{s}(\bar{ \varphi } )=0$.
So we have to decide if to include formulas in which this occurs.
Does 
\end{remark}

\begin{PROOF}{\ref{b13}}
Consider a sentence $\psi \in \bbL(\mathbf Q)$, see \ref{b8}.
\mn
\begin{enumerate}
\item[$\boxplus_0$]  for each $n$ we consider drawing 
$(\cG_{p,n},\bar{\mathbf t}) \in \text{ Graph}_n \times \mathbf T$, that
is, independently we draw
\sn
\begin{enumerate}
\item[$\bullet$]   $\bar{\mathbf t} \in \mathbf T$ by the probability space from
  \ref{b2}(2)
\sn
\item[$\bullet$]  $\cG_{n,\bar p} \in \Graph_n$ = the set of graphs with set of
  nodes $[n]$ with each edge drawn with probability $p_n$
independently of the other edges
\end{enumerate}
\sn
\item[$\boxplus_1$]  It suffices to prove that
\sn
\begin{enumerate}
\item[$(a)$]   the probability of $``\cG_{n,p} \models_{\bar{\mathbf t}}
  \psi"$, i.e. the pair $(\cG_{n,\bar p},\bar{\mathbf t})$ satisfies
this, either is $\ge \frac{1}{2^{nr}}$ or is
$\ge 1 - \frac{1}{2^{nr}}$ for some $r = r(\psi) \in (0,1)_{\bbR}$
\sn
\item[$(b)$]   which case does not depend on $n$
\sn
\item[$(c)$]   moreover the probability is $\ge 1 -
  \frac{1}{2^{n\ge}}$ iff $\psi \in T_*$.
\end{enumerate}
\end{enumerate}
\mn
[Why?  Consider the drawing of $(\langle \cG_{n,p}:n \in
\bbN\rangle,\mathbf t) \in \prod\limits_{n} \Graph_ n
    \times \mathbf T$.
For every $\psi \in \bbL(\mathbf Q)$, the following event
$\cE^1_\psi \wedge \cE^2_\psi$ has probability zero, where

\[
\cE^1_\psi := (\text{for infinitely many } n,\cG_{n,p} \models_{\mathbf
  t} \psi)
\]

\[
\cE^2_\psi := (\text{for infinitely many } n,\cG_{n,p} \models_{\mathbf
  t} \neg \psi).
\]

\mn
This holds by (a)+(b) of $\boxplus_1$.  Hence also the event $\cE =
\bigvee\{\cE^1_\psi \wedge \cE^2_\psi:\psi \in \bbL(\mathbf Q)\}$ has
probability zero.  Hence, by Fubini theorem, drawing for a set of
$\bar{\mathbf t}$'s of measure 1, the event $\cE^1_\varphi[\mathbf t]
\wedge \cE^2_\psi[\mathbf t]$ has probability zero, where
$\cE^\ell_\psi[\mathbf t]$ is the event $\cE^\ell_\psi$ fixing $\mathbf t$.]

To prove $\boxplus_1$, fix $\psi \in \bbL(\mathbf Q)(\tau_{\gr})$.
We can find a $\bar\Delta$ such that:
\mn
\begin{enumerate}
\item[$\boxplus_2$]
\begin{enumerate}
\item[(a)]  $\bar\Delta = \langle \Delta_\ell:\ell
  \le \ell(*)\rangle$
\sn
\item[(b)]  $\Delta_\ell$ is a finite set of formulas from
  $\bbL(\mathbf Q)$ increasing with $\ell$
\sn
\item[(c)]  $\Delta_0$ is the set of quantifier free
  formulas
\sn
\item[(d)]  $\psi \in \Delta_{\ell(*)}$
\sn
\item[(e)]   every formula in $\Delta_{2 \ell +1}
  \backslash \Delta_{2 \ell}$ is
  gotten from formulas from 
$\Delta_{2 \ell}$ by a first order operation
  $(\neg \varphi(\bar x),\varphi_1(\bar x) \wedge \varphi_2(\bar x),
\exists y \varphi(\bar x,y))$
\sn
\item[(f)]   every formula in $\Delta_{2 \ell +2}
\backslash \Delta_{2 \ell +1}$ is of the form $\psi(\bar z) = (\mathbf Q
\ldots \bar x_i,\bar x'_i,\ldots)_{i<k} \bar\varphi(\bar z)$
where $\bar\varphi = \bar\varphi(\bar z)$ recalling
\ref{a51} is a complete reduced $\mathbf s'$-scheme 
for some $\mathbf s'$, i.e. is as in Definition \ref{a6} but the
$\varphi_{0,i}(\bar x,\bar z),
\varphi_{1,i,j}(\bar x_i,\bar x_j,\bar z),\varphi_2(\bar z)$ 
being from $\Delta_{2 \ell +1}$
\sn
\item[(g)]   no two distinct $\bar\varphi$'s which occur in
  $\Delta$ on $(\mathbf Q,\ldots)\bar\varphi$ are explicitly isomorphic
(see Definition \ref{a45}), \but \, replacing equality of formuals by
  equivalence for every randm enough $\cG_{p,n}$ (during the proof
  this will get a syntactical characterization).
\end{enumerate}
\end{enumerate}
\mn
[Why?  Should be clear.]
\mn
\begin{enumerate}
\item[$\boxplus_3$]  let $\Delta_\ell = \{\vartheta_s(\bar x_s):s \in
  I^*_\ell\}$ hence $\cT_\ell$ is finite and
 $m < \ell \Rightarrow I^*_m \subseteq I^*_\ell$.
\end{enumerate}
\mn
Now by induction on $\ell \le \ell(*)$ we choose $\mathbf s_\ell,
\bar\vartheta'_\ell,\bar\vartheta''_\ell$ and the function 
$G \mapsto M_{G,\bar{\mathbf t},\ell}$ for $G$ a graph on $[n]$ some $n$ 
such that:
\mn
\begin{enumerate}
\item[$\boxplus_{4,\ell}(A)$]  
\begin{enumerate}
\item[(a)]  $I_\ell$ finite
\sn
\item[(b)]  $\mathbf s_\ell$ is as in Definition
  \ref{a2}, an $I_\ell$-kind of a vocabulary
\sn
\item[(c)]
\begin{enumerate}
\item[$(\alpha)$]  $\mathbf s_0,I_0$ are defined by $I_0 =
  \{s_0\}$ for some $s_0 \notin I^*_{\ell(*)},
n_{\mathbf s,s_0} = 2,K_{\mathbf s_0,s_0}
  =\Sym(2)$, the group of permutation of $\{0,1\}$
\sn
\item[$(\beta)$] $I_{2 \ell +1} = I_{2 \ell}$
\sn
\item[$(\gamma)$]  $I_{2 \ell+2} =I_{2 \ell +1} \cup
  (I^*_{2 \ell +2} \backslash I^*_{2 \ell+1})$
\sn
\item[$(\delta)$]   so
\begin{itemize}
\item  $\langle I_\ell:\ell \le \ell(*)\rangle$ is increasing
\sn
\item  $\mathbf M_{\mathbf s_0,n}$ is
  $\Graph_n$, the set of graphs with set of nodes $[n]$
\end{itemize}
\end{enumerate}
\sn
\item[(d)]  $\bar\vartheta'_\ell = \langle
  \vartheta'_s(\bar x_s):s \in I_\ell\rangle$
\sn
\item[(e)]  $\vartheta'_s(\bar x_s)$ a formula in
  $\bbL(\tau_{\mathbf s_\ell})$ for $s \in I_\ell$
\sn
\item[(f)]  $\vartheta''_s(\bar x_s)$ is a quantifier free
  formula in $\bbL(\tau_{\mathbf s_\ell})$ equivalent to
  $\vartheta'_s(\bar x_s)$ in 
 the limit theory $T_{\mathbf s_\ell}$, see Definition \ref{a4}
\sn
\item[(g)]   for any given $G \in \cG_{n,p}$, i.e. $G \in
  \mathbf M_{\mathbf s_0,n}$ and $\bar{\mathbf t} \in \mathbf T$ 
we define
 $M_{G,\bar{\mathbf t},\ell} \in \mathbf M_{\mathbf s_\ell,n}$ by:
\sn
\begin{enumerate}
\item[$(\alpha)$]  $M_{G,\bar{\mathbf t},\ell}$ is a
  $\tau_{\mathbf s_\ell}$-model expanding $M_{G,\bar{\mathbf t},m}$ for $m
  < \ell$ \underline{and} for 
$s \in I_\ell,R^{M_{\cG,\bar{\mathbf t},\ell}}_s$ is 
defined by $\vartheta_s(\bar z_s)$ and also by $\vartheta'_s(\bar z_s)$
\sn
\item[$(\beta)$]   if $\ell=0,M_{G,\bar{\mathbf t},\ell}$ is $G$
\sn
\item[$(\gamma)$]   if $\ell=2m+1,s \in I_\ell
  \backslash I_{2m}$ we apply the the first order 
construction of $\vartheta_s(\bar x_s)$ from the 
formulas $\langle \vartheta_s(\bar x_s):s \in I_{2m}\rangle$
to construct $\vartheta_s(\bar z_s)$ from
  $\langle \vartheta''_t(\bar x_t,\bar z_s):t \in I_{2m}\rangle$
\sn
\item[$(\delta)$]   for $\ell = 2m+2$ and $s \in I_\ell
  \backslash I_{2m+1}$ if $\vartheta_s(\bar x_s) = (\mathbf Q \ldots,\bar
  x_i,\bar x'_i,\ldots)_{i < \mathbf i(\bar\varphi_{\mathbf i})}
  \bar\varphi_s(\bar x_s)$ 
we define $\vartheta'_s(\bar x_s)$ by
  replacing in $\bar\varphi_s$ every $\vartheta_t$ by $\vartheta''_t$
  getting $\bar\varphi'_s$ 
and let $\vartheta'_s(\bar x_\iota) =
  (\mathbf Q \ldots,\bar x_i,\bar x'_i,\ldots)_{i < \mathbf
  i(\bar\varphi_\iota)} \bar\varphi'_s$
\sn
\item[$(\varepsilon)$]   we choose $\vartheta''_i(\bar
  x_i)$ by \ref{a4} sequence clause (f) here.
\end{enumerate}
\end{enumerate}
\end{enumerate}
\mn
Now for each $\ell \le \ell(*)$ we have two relevant ways to draw as
$\mathbf s_\ell$-structure $M$ withuniverse = set of elements $[n]$.

First, draw $\mathbf t \in \mathbf T$ and $\cG = \cG_{q,n}$ (recall $q \in
(0,1)_{\bbR}$ was fixed in the beginning of Theorem \ref{b13}) and
compute $M_{\cG,\bar{\mathbf t},n}$, a $\mathbf s_\ell$-structure.  This
induces a distribution $\mu_{q,n,\ell}$ on $\mathbf M_{\mathbf s_\ell,n}$,
i.e. $\mu_{q,n,\ell}(M) = \Prob(M_{\cG_{q,n},\mathbf t,n} =
M|\mu_{\gr,q,n} \times \mu_{\mathbf T} = \mu_{\mathbf s_{\gr},\bar
  p_{\gr,q,n}} \times \mu_{\mathbf T})$.

Second, we shall choose $\bar p_\ell \in \mathbf P^2_{\mathbf s_\ell}$ and
draw $\cM_{\mathbf s,\bar p_\ell,n}$ here the distribution is ?  The
interest in the first is that our aim is to prove the 0-1 law for
$M_{\cG,\bar p,n}$, in particular, for $\ell=\ell(*)$ and our sentence
$\psi$; we use the other $\ell$'s in an induction.

A priori the probability of $``M_{\cG,p,n} \models \psi"$ is opaque.

For the second, $\cM_{\mathbf s,\bar p_\ell,n}$ an understanding of the
probability of $M_{\mathbf s,\bar p_\ell,n} \models \psi$ is now well
known and satisfies the 0-1 law.  Hence it suffices to prove that the
distribution of $M_{\cG,\mathbf t,\ell}$ (for $\cG \in \cG_{p,n})$ from
$\mathbf M_{\mathbf s_\ell,n}$ and $M_{\mathbf s_\ell,\bar p,n} \in \mathbf
G_{\mathbf s_\ell,n}$ are sufficiently similar.

Naturally we choose:
\mn
\begin{enumerate}
\item[$(*)_1$]  $(a) \quad p_{\mathbf s_{\gr},s_{\gr},n} = p_{\mathbf
    s_0,s_{\gr},n} = q$
\sn
\item[${{}}$]  $(b) \quad p_{\mathbf s_\ell,t,n} = q/\mathbf g(n)$ for 
    $ t \in I_{\mathbf s_\ell} \backslash \{s_{\sg}\}$.  
\end{enumerate}
\mn
Of course, we induct;  
    for $\ell=0$ there is no difference so we
deal now with $\ell+1$ if $\ell$ is even this is trivial so assume
$\ell$ is odd.

There are several reasons for a difference, for a given model $M \in
\mathbf M_{\mathbf s_\ell,n}$
\mn
\begin{enumerate}
\item[$(*)_{M,1}$]  $t \in I^*_{\ell +1} \backslash I^*_\ell$ and
$\bar c \in \varphi_{t,2}(M)$.  The graph $H_{\bar\varphi_t,M,\bar c}$ 
is $\iota$-low (for a given $n$ there are at most
$n^{k(\bar\varphi'_t)}$ (check cases)
\sn
\item[$(*)_{M,2}$]  for some $t(1),t(2) \in I^*_{\ell +1} \backslash
  I^*_\ell,\bar c_2 \in \varphi_{t(1),2}(M)$ and $\bar c_2 \in
  \varphi_{t,2}(M)$ we have $(t(j),\bar
  c_1/E_{\bar\varphi'_{t(1)}} \ne (t(2),\bar
  c_2/E_{\bar\varphi'_{t(2)}})$ but the graphs
  $H_{\bar\varphi'_{t(1)},M,\bar c_1},H_{\bar\varphi'_{t(2)},M,\bar
    c_2}$ are isomorphic
\sn
\item[$(*)_{M,3}$]  for some $t(1),t(2) \in I^*_{\ell +2} \backslash
I^*_\ell$ and $t(2) \in \cup\{I^*_{2k+2} \backslash I^*_{2k+1}:2k+2
\le \ell\}$ and $\bar c_1 \in \varphi_{t(1),2}(M),\bar c_2 \in 
\varphi'_{t(2),2}(M)$ the graphs 
$H_{\bar\varphi'_{t(1)},M,\bar c_2},H_{\bar\varphi'_{t(2)},M,\bar
  c_2}$ are isomorphic
\sn
\item[$(*)_{M,4}$]  the sequence $\bar p \in \mathbf P^2_{\mathbf q}$ try
  to immitate $\mathbf t$, \underline{but} having the probability for
  $\cM_{\mathbf s_{\ell +1},\bar p,n} \models R_t[\bar c]$ is $p_{t,n} =
  1/\mathbf g(n)$ whereas the probability $\mathbf t_i =1$ is $1/\mathbf
  g(|H^*_i|)$ where $i$ is such that
$H_{\bar\varphi_t,\cM_{\cG,\mathbf t},\bar c} = H^*_i$ for $\cG =
{\mathscr G} _{q,n}$.   
\end{enumerate}
\mn
Now there is no reason that usually $i=n$.  However, if $\iota = 2$
then $|H^*_i| \le k(\bar\varphi_t) \cdot n$ and if $\iota=1,H^*_1 \le
n^{k(\bar\varphi_2)}$.  In both cases with probability very close to
1, (for $\mu_{\mathbf s_{\ell +1},\bar p,n}$), $|H^*_i| \ge
n/2^{k(\bar\varphi_t)}$.  So clearly as $\mathbf q$ grow slowly enough,
see \ref{a0}(2).

This is also true for $(*)_{M,1},(*)_{M,2},(*)_{M,3}$.  Together, we
have two distributions on $\mathbf M_{\mathbf s_{\ell+1},n}$ and for the
second, omitting a set of $M$ with small probability (in $\mu_{\mathbf
  s_{\ell +1},\bar p,n}$) for any other $M$, the two distributions
give almost the same values.  The computations are easy so we are done.
\end{PROOF}

\begin{remark}
\label{b20}
To eliminate $(*)_4{M,}$   
        in the end of the proof we may complicate the
drawing of $\cM_{\mathbf s_{\ell+1},\bar p,n}$  
We draw $\cM_{\mathbf
  s_m,\bar p,n}$ by induction on $m$: if $m=2j+2,M = M_{\mathbf
  x_{2j+1},\bar p,n}$ given for $R_t(t \in I^*_m \backslash
I^*_{2k+1})$ we consider only $\bar c \in \varphi'_{t,2}(M)$ let $m =
m_t(\bar c)= m_t(\bar c,M)$ be the number of nodes of
$H_{\bar\varphi'_t,M,t}$ and we draw a truth value of $R_t(\bar c)$
with probability $1/g(m)$.  Proving the 0-1 law for such drawing is easy.
\end{remark}
\newpage

\section {How to get a real quantifier, i.e. definable $K$}

\begin{discussion}  \label{d1}
In the introduction we have considered drawing a truth value to all graphs.
So replacing \lqq converge to zero or to one" we ask only
\lqq for every $ \varepsilon > 0 $ for every $ n $ large enough
the probability is up to $ \varepsilon $ closed to zero or to one,
The point is that otherwise we can weakly express
``$|\varphi_1(\cG_{p, n } 
,\bar a_1)| = |\varphi_2(\cG_{p,n},\bar a_2)|$,
e.g. for $\varphi(x,y) = xRy$.  So we can find $\psi_1(x_1,x_2)$
implying $\valency_{\cG_{p,n}}(y_1) = \valency_{\cG_{p,n}}(y_n)$, this
will complicate the matter.

In more details, let $\psi_\varphi  
    (y)$ say ``the empty graph on
$\varphi(\cG_{p,n},y)$ is green".   

Let $\psi_2(y_1,y_2)$ say:
\mn
\begin{enumerate}
\item[$(a)$]  $\psi_1(y_1) \equiv \psi_2(y_2)$
\sn
\item[$(b)$]  for $\ell \in \{1,2\}$ and $y'_\ell$ there is
  $y_{3-\ell}$ such that $|\varphi(\cG_{p,n},y_1) \cap
 \varphi(\cG_{n,p},y'_1)| = |\varphi(\cG_{p,n},y_2) \cap
  \varphi(\cG_{p,n},y'_2)|$.
\end{enumerate}
\mn
This nearly expresses $|\varphi(\cG_{p,n},y_1)| =
|\varphi(\cG_{p,n},y_2)|$.  We can strengthen this and find
approximation to $a+1$ and cases of addition.

While the above does not suffice to prove impossibility, it suffices
to show the problem is not promising and is different; 
    maybe relevant is the late \cite{Sh:F1943}.  
    
\end{discussion}
    
\begin{discussion} \label{d4}
Can we use a quantifier $\mathbf Q_{\mathbf K}$
which depends just on the number of edges via the number of nodes.

\noindent
1) If it depends only on the number of nodes, 
it seemed 
that this is bad for 0-1 laws.

\noindent
2) Notes that surely graphs $H_1,H_2$ occur
up to isomorphism 
when
$H_2$ is gotten by omitting one edge of $H_1$.  So we may try that it
depends  
    only 
the number modulo $(\lfloor \log \log(4+1)\rfloor)!$
Quite reasonable choice of the quantifier 
but not ideal.  

\noindent
3) So we may try to change the logic such that 
essentially just 
changing one edge does
not matter; that is excluding some family of graphs which 
with probability  
one does not occurs for a random enough
$ {\mathscr G} _{p,n}$.  
This is a reasonable logic, even without ``$H \in \mathbf
K$ depends just on the number of edges (and nodes)"
\mn
\begin{enumerate}
\item[(A)]  if we forget this restriction, we need to change the
  flipping of coins for the logic, e.g. fixing size first, choose one
  randomly, do this for each neighborhood, choose with distorted
  probability; not clear if converge and 
  there is a natural way
\sn
\item[(B)]  Here  
$n^{\mathbf g(h)}$ goes slowly to $\infty$ and is used
  how to make the results O.K..  
  Note: in ${\mathscr G} 
  _n$ the size of a
  definable graph for some $m,$ is $\approx \frac{n}{m}$ so the
  variance is $c \sqrt\frac{n}{m}$; still the edges have probability
  $\frac 12$ and so O.K.
\end{enumerate}
\mn
However  
for later $M^\iota_n$ ($\iota <$ quantifier depth) the probability
of each case of a relation is, i.e. $H \in \mathbf K$ for a structure
with probability $\frac{1}{\mathbf h(n)}$ so manipulating $\mathbf h$
gives different results.

\noindent
4) But we have a more profound problem: we have nicely definable
$H_1,H_2$ 
getting $ H_2 $ from $  H_1$  by, for some nodes 
$a \not= b$ 
omitting the edges $ (a, c) $ and adding the 
edges  
$ (b, c)$ whenever $ (b, c))$ is an edge  when 
$ (a,c)$ is not.



Alternatively 
omitting the edge $(a,c)$ when $(b,c)$ is an edge, 
The first does not change  the 
number of edges, the second changes seriously.  This may be close to
the variance for the number of edges.

A medicine? ask: omitting $\log_*(H_)$ edges, what is the minimal
number of edges?

The overcoming may cost: in how to make the probability 
computations 
right.

\noindent
5) Note: from random ${\mathscr G} 
_{n,1/2}$ we build $ {\mathscr M } 
^9_{n_1} = {\mathscr G} 
_{n,H}$ an 
$\mathbf s_0$-structure $
{\mathscr M }  
^1_n$ expands by applying the quantifier
getting an $\mathbf s_1$-structure.  But $
  {\mathscr M }  
  ^{\mathbf s_1}_n$ is different:
\mn
\begin{enumerate}  
\item[(a)]  for ${\mathscr M }  
^{\mathbf s_1}_n$ the cases are totally independent
\sn
\item[(b)]  ${\mathscr M } 
^1_n$ is different: first we draw $R_{\gr}$ ($= R_0$ in
  the lecture) after this we draw the other relatives but their
  probabilities:
\sn
\begin{itemize}
\item  depends on the drawing of 
${\mathscr M } _n = {\mathscr G} _{1/2, n}$  
\sn
\item  in particular, on the sizes of the $H$'s which are not too far
  from $n$ but are different.
\end{itemize}
\end{enumerate}
\mn
This complicates our work but the estimates are not so different.
\end{discussion}

\begin{discussion} \label{d7}
One which seems easiest while not unreasonable is: given a finite
graph $G$, with $m$ points, which is reasonable - defined as in
\cite{Sh:F1166} and a point $b$ in it, compute the valency minus
$m/2$, divided by square root of $m$ (or the variance of the related
normal distribution) and ask if rounding to integers is odd or even.

 We may replace the valency by the number of edges of $G$.

What are the dangers?  As we may define a variant of the graph
omitting one edge, in some cases this will change the truth value.
For each nod the probability goes to zero but in binomial distribution
the probability of e.g. getting valency exactly half of the expected
value (rounded) is about 1 divided by the square root of $m$.

So we should divide not by the square root of $m$ but by a larger
value (maybe instead of asking on even/odd of the rounded value just
ask if it can be larger than one, or absolute value) such that:
\mn
\begin{enumerate}
\item[$(a)$]  almost surely (i.e. with large probability) for some
  node the value is above 1
\sn
\item[$(b)$]  the probability that it is exactly one for some node is
  negligible, and this is true even if we use a graph only definable
  (reversing edge/non-edge, omitting some, etc.).
\end{enumerate}
\mn
So we should say that clearly by continuity considerations there are
such choices.  A danger is that the $n$ being odd/even can be
expressed.

Another avenue is to choose the more natural ``the valency is at least
half"; but then it seems we can express being even/odd: say change by
one edge change the truth value and this is true even if we omit one
node.  So the number of neighborhoods is half in both cases.
\end{discussion}
\newpage


\bibliographystyle{amsalpha}
\bibliography{shlhetal}

\end{document}